\documentclass[12pt,a4paper,reqno]{amsart}

\usepackage{amssymb}
\usepackage{amscd}
\usepackage{amsfonts}
\usepackage{graphicx}
\usepackage{asymptote}
\usepackage{hyperref}
\usepackage[british]{babel}
\usepackage{enumerate}

\numberwithin{equation}{section}

\theoremstyle{plain}
  \newtheorem{theorem}{Theorem}[section]
  
  \newtheorem{lemma}[theorem]{Lemma}

\theoremstyle{remark}

\renewcommand{\leq}{\leqslant}
\renewcommand{\geq}{\geqslant}

\newcommand\PP{\mathbb{P}}
\newcommand\QQ{\mathbb{Q}}
\newcommand\E{\mathbb{E}}
\newcommand\Z{\mathbb{Z}}
\newcommand\R{\mathbb{R}}

\newcommand\C{\mathbb{C}}

\newcommand\N{\mathbb{N}}
\newcommand\W{\overline{W}}

\newcommand\eps{\varepsilon}
\newcommand\vol{\operatorname{vol}}

\newcommand\w{\mathrm{w}}

\newcommand{\Mod}[1]{\;(\mathrm{mod}\,#1)}

\renewcommand\>{\rangle}

\begin{document}
\title[Correlations of representation functions of binary quadratic forms]
{Correlations of representation functions \\ of binary quadratic forms}

\author{Lilian Matthiesen}
\address{School of Mathematics\\
University Walk\\
Bristol, BS8 1TW\\
United Kingdom
}
\email{l.matthiesen@bristol.ac.uk}
\thanks{While working on this paper the author was supported by EPSRC
grant number \texttt{EP/E053262/1}.}

\begin{abstract}
The purpose of this short note is to extend previous work on
linear correlations of representation functions of positive definite
binary quadratic forms to allow indefinite forms.
\end{abstract}

\maketitle
\tableofcontents

\addtocontents{toc}{\protect\setcounter{tocdepth}{1}}

\section{Introduction}
This paper should be thought of as an appendix to \cite{m-definiteforms}.
Let $f: \Z^2 \to \Z$ be a primitive binary quadratic form whose
discriminant is not a square.
The latter assumption is equivalent to the assertion that $f$ does not
factor over $\QQ$.
Thus $f$ is irreducible and there is no non-trivial representation of $0$.
We shall be concerned with the function $r_f: \Z \to \N_0$ that counts
for each integer $n$ the number of primary representations of $n$ by $f$.
If $\mathcal E_f$ denotes the group of automorphs of $f$, then the
function $r_f$ takes at $n\in\Z$ the form
$$r_f(n) = |\{ (x,y) \in \Z^2/\mathcal{E}_f: f(x,y) = n \}|~.$$
We shall assume familiarity with the notation from \cite{m-definiteforms}.
Extending results from \cite{m-definiteforms}, our aim is to establish the
following theorem.

\begin{theorem}\label{thm:main}
Let $f_1, \dots, f_t$ be primitive irreducible binary quadratic forms.
We further assume that every definite form is positive definite.
For each $f_i$ let $r_{f_i}$ denote the representation function defined
above, and let $D_i$ be the discriminant of $f_i$.
If $D_i<0$, let $\w(D_i)= 6,4,2$ according as $D_i=-3,-4$ or $D_i<-4$.
If $D_i>0$, let $\eps_i$ denote the fundamental unit of $\QQ(\sqrt{D_i})$.

Let $\Psi= (\psi_1, \dots ,\psi_t): \Z^d \to \Z^t$ be a system of
affine linear forms such that no two forms $\psi_i$ and $\psi_j$ are
affinely dependent, and suppose that $K \subset [-N,N]^d \subset \R^d$ is
a convex body such that $\psi_i(K) \geq 0$ whenever $f_i$ is
positive definite.
Then 
\begin{align*}
 \sum_{n\in \Z^d\cap K } r_{f_1}(\psi_1(n)) \dots r_{f_t}(\psi_t(n)) 
= \beta_{\infty} 
   \prod_{p} \beta_p
 +  o(N^d)~,
\end{align*}
where the implicit constant may depend on $D_1,\dots,D_t$ and the
coefficients of $\Psi$.
The constants are given by
$$\beta_{\infty}
= \vol(K)
  \prod_{i: D_{i} < 0} 
  \frac{2\pi}{\w(D_i)\sqrt{-D_i}}
  \prod_{j: D_j > 1} 
  \frac{\log \eps_j}{\sqrt{D_j}}~,
$$
and
$$\beta_p 
= \lim_{m \to \infty}
  \E_{a \in(\Z/p^{m}\Z)^{d}} 
  \prod_{i = 1}^t
  \frac{\rho_{f_i,\psi_i(a)}(p^{m})}{p^{m}}~,$$
with $\rho_{f,A}(q)$ denoting the local number of representations of
$A \Mod{q}$ by $f$, that is,
$$\rho_{f,A}(q)
:= |\{(x,y)\in (\Z/q\Z)^2:f(x,y) \equiv A \Mod{q}\}|~.$$ 
\end{theorem}
Theorem \ref{thm:main} provides a useful tool in the study of several
Diophantine problems.
It will prove instrumental, for example, in forthcoming joint work of the
author with Browning and Skorobogatov on the Hasse principle and weak
approximation for conic bundle surfaces.
Generalising work of Heath-Brown \cite[Theorem 2]{heath-brown}, it can
also be used to study the set of $\QQ$-rational points on varieties 
$X \subset \PP^{2r+1}$ of the form
$$Q_i(Y_i,Z_i)=L_i(U,V)M_i(U,V), \quad 1 \leq i \leq r.$$
Here each $Q_i$ is an irreducible binary quadratic form over $\QQ$ and
$L_i,M_i$ form a set of pairwise non-proportional linear forms defined
over $\QQ$.
Mimicking the approach of Heath-Brown, which dealt with the case $r=2$
when $Q_i(Y,Z) = Y^2+Z^2$ for $i=1,2$,
it would be straightforward to deduce from Theorem \ref{thm:main} that
the set $X(\QQ)$ is Zariski dense in $X$ as soon as there exists a
non-singular $\QQ$-rational point.
It is worth stressing that very little is known about the arithmetic of
intersections of $r\geq3$ quadratics when the ambient projective
dimension $n$ is small.
A notable exception is found in work of Skorobogatov \cite{sk} who
handles the question of weak approximation when $r=3$ and $n\geq11$,
under suitable hypotheses. 

Following the strategy from \cite{m-divisorfunction, m-definiteforms}, we
deduce Theorem \ref{thm:main} by means of the machinery Green and Tao
developed in 
\cite{
green-tao-linearprimes,
green-tao-nilmobius,
green-tao-polynomialorbits},
in connection with Green, Tao and Ziegler's inverse theorem for the
uniformity norms \cite{gtz}.
In fact, it will be possible to modify the approach from
\cite{m-definiteforms}, adapting only few parts.
We therefore restrict attention to these changes, referring the reader to
\cite{m-definiteforms} for details regarding the remaining parts.

Our main tasks in establishing Theorem \ref{thm:main} are
then to construct a pseudorandom majorant for $r_f$ when $f$ is
indefinite, to show that a suitably $W$-tricked version $r'_f$ of this
function is Gowers-uniform, and to check that the linear forms and
correlation estimates are valid across the different majorants from the
definite and indefinite cases.

Theorem \ref{thm:main} can be extended to reducible forms.
If $f$ is a form of discriminant $1$, then the representation function of
$f$ is given by
$\frac{1}{2}\tau(n) = \frac{1}{2} \sum_{d \in \N} 1_{d|n}$ for non-zero
integers $n$. 
If the discriminant is a square greater than $1$, then the representation
function is a restricted divisor function.
In both cases, the majorant function from \cite{m-divisorfunction} can be
employed.
The $W$-trick is compatible with the one used for irreducible forms.
The proof of the linear forms condition can be adapted without
difficulties (although this is a tedious task and will not be treated
here).

All implicit constants in this paper are allowed to depend on the
discriminants $D_1,\dots,D_t$ and the coefficients of $\Psi$.

\section{Representation by indefinite forms}
Let $f(x,y)=ax^2+bxy+cy^2=\<a,b,c\>$ be a primitive indefinite form whose
discriminant is not a square.
In order to obtain a more explicit expression for $r_f$ than that given
in the introduction, we proceed to fix a fundamental domain
$K_0 \subset \R^2$ for the action of the units in the associated real
quadratic field. This will allow us to write
$$r_f(n) = |\{ (x,y) \in K_0: f(x,y) = n \}|~.$$
 
Let $D=b^2 -4ac$ denote the discriminant of $f$ and 
let $\eps = \frac{1}{2}(t_0 + u_0 \sqrt{D})$ be the fundamental unit of
$\QQ(\sqrt{D})$, that is, $(t_0,u_0)$ is the solution of $t^2 - Du^2 = 4$
for which $u_0$ is minimal under the condition $t_0,u_0>0$.
Then, cf. \cite[\S9.3]{rose}, the automorph 
$$T_0 := \frac{1}{2}
\begin{pmatrix}
 t_0-bu_0 & 2au_0 \cr
 -2cu_0   & t_0+bu_0
\end{pmatrix}
$$
generates the group $\mathcal{E}_f = \{\pm T_0^n \mid n \in \Z \}$.
In what follows we restrict attention to positive $n$; on
replacing $f$ by $-f$, this also covers the case of representation of
negative integers by $f$.
Since equivalent forms have identical representation functions, we
may assume that $f=\<a,b,c\>$ satisfies
$|b|\leq |a|,|c| \leq \frac{1}{2}\sqrt{D}$ and
$a>0>c$.
Thus, either $\<a,b,c\>$ or $\<c,b,a\>$ is reduced.
The automorph $T_0$ takes the line $\{y=0\}$ to $\{y= \theta x\}$ for
$\theta = -2cu_0/(t_0-bu_0)$.
Observe that $\theta$ is bounded in terms of $D$, and satisfies $\theta>0$
since $c<0$ and $t_0-bu_0 > t_0-\sqrt{D}u_0 > 0$.
Let $$K_0 \subset \{(x,y) \in \R^2: x > 0, y \geq 0 \}$$  be the cone
defined by the lines $\{y=0\}$ and $\{y= \theta x\}$. Then
$$r_f(n) = |\{(x,y) \in K_0 \cap \N^2: f(x,y)=n\}|$$
for $n>0$.

If $K_0(N):=\{(x,y) \in K_0: f(x,y) \leq N\}$, then 
(cf.~\cite[\S10, Lemma 3.5]{rose})
$$\vol K_0(N) = N \frac{\log \eps}{\sqrt{D}}~.$$

Thus, the analogue of \cite[Lemma 4.1]{m-definiteforms} is the following.
\begin{lemma}\label{lem:average-in-APs}
 Let $q>0$ and $A$ be integers and let 
$P:= \{n\leq N: n \equiv A \Mod{q}\}$ 
be an arithmetic progression.
Then the average of $r_f$ along $P$ satisfies
$$\E_{n\in P} r_f(n) 
= \frac{\log \eps}{\sqrt{D}}
  \frac{\rho_{f,A}(q)}{q}
  + O(|P|^{-1/2} q^2)~.$$
\end{lemma}

\section{$W$-trick and Pseudorandom majorant}

As in the definite case, $r_f$ has the arithmetic representation
$$r_f(n)= \sum_{d|n} \chi_D(d)$$
for non-zero $n$ which are coprime to $D$.
All non-archimedean results that are recorded in 
\cite[\S5--6]{m-definiteforms} carry over to indefinite forms directly, as
is to be expected since definiteness and indefiniteness are archimedean
properties.
In particular, 
$$r_f(n) \ll \sum_{d|n} \chi_D(d)~.$$
holds for any non-zero integer $n$.
This allows us to employ a majorant of the exact same form as in the
definite case, see \cite[\S2, \S7.2]{m-definiteforms}.

Let $C_1>1$ be the parameter that appears in the divisor function majorant
\cite[Proposition 2.3]{m-definiteforms} and assume it to be sufficiently
large for the conclusion of \cite[Lemma 3.3]{m-definiteforms} to hold.
We then define, as in the definite case, $\W := \prod_{p<w(N)}
p^{\alpha(p)}$,
where $w(N)=\log\log N$ and
$p^{\alpha(p)-1} < \log^{C_1+1} N \leq p^{\alpha(p)}$.
For integers $0 \leq A < \W$, write $r_{f,A}(n) := r_f(\W n + A)$.
If $\rho_{f,A}(\W)>0$, we may in view of Lemma \ref{lem:average-in-APs}
define the normalised function $r'_{f,A}:\N \to \R_{\geq 0}$ by
$$r'_{f,A}(n) := 
\Big(\frac{\log \eps}{\sqrt{D}} \frac{\rho_{f,A}(\W)}{\W}\Big)^{-1}
r_{f,A}(n).$$
Since the major arc estimate \cite[Proposition 7.4]{m-definiteforms} only
builds on non-archimedean properties, it continues to hold for indefinite
forms $f$:
\begin{lemma}\label{lem:major_arc}
Let $P \subseteq [N/\W]$ be a progression of $w(N)$-smooth common
difference $q_1$ and let $0<A<\W$ be such that $\rho_{f,A}(\W)>0$ and
$A\not\equiv0\Mod{p^{\alpha(p)}}$ for all primes $p<w(N)$. 
If $P = \{q_1 m + q_0: 0 \leq m < M \}$ has length $M$, then
$$
  \E_{n \in P} r'_{f,A}(n) 
= \E_{0 \leq m < M} r'_{f,A}(q_1m+q_0)
= 1 + O\Big(\frac{\W(\W q_1)^2}{M^{1/2}}\Big)~.
$$
\end{lemma}
Lemma \ref{lem:major_arc} shows that the same $W$-trick works for
definite and indefinite forms and reduces matters to studying the
functions $r'_{f,A}$ and their majorants $\nu'_{D,\gamma}(\W n +A)$
defined in \cite[\S7.2]{m-definiteforms}.

If all forms in the collection $f_1,\dots,f_t$ are irreducible,
then the required linear forms and correlation estimates for the
corresponding $\W$-tricked versions of the majorants follow
word by word the proofs in \cite[\S9]{m-definiteforms}. 
Some minor changes would be required if reducible forms are included.
 
\section{The non-correlation estimate}
Let $h:\N \to \C$ be an arithmetic function.
In order to prove that
$$ \E_{n\leq N} h(n) F(g(n)\Gamma) = o_{G/\Gamma}(1)$$
for all nilmanifolds $G/\Gamma$ and all polynomial nilsequences
$(F(g(n)\Gamma))_{n \in \N}$ on it, it usually suffices to
show that such an estimate holds when $\int_{G/\Gamma} F = 0$ and when
the nilsequence is equidistributed in some quantitative sense.
This strategy was introduced by Green and Tao in
\cite{green-tao-nilmobius} and employed in \cite{m-definiteforms},
although several modifactions were necessary in the latter case. 
The strategy that was used for $h(n)=r'_{f,A}(n)$ in the case of
definite forms $f$ continues to apply in the indefinite case with the only
changes concerning the final part of the argument carried out in
\cite[\S18]{m-definiteforms}.
There the summation  
$$
\E_{n\leq N'} r_f(\W n + b) F(g(n)\Gamma)
= \frac{1}{N'} 
  \sum_{\substack{x,y \geq 0: \, f(x,y) \leq N, \\
                  f(x,y) \equiv b \Mod{\W} }}
  F(g((f(x,y) - b)/\W)\Gamma),
$$
for $N':= [(N-n)/\W]$, is split into ranges where either $x$ or $y$ is
fixed, while the free variable varies over a long enough interval for an
analogue of Weyl's inequality to apply.

\begin{figure}
\centering
\includegraphics{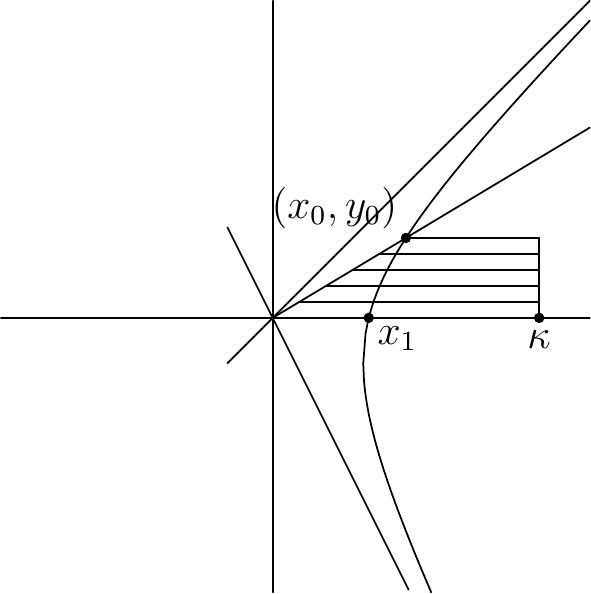}

\caption{The summation is split into individual sums for fixed $y$. Each
such sum is expressed as the difference of two sums over a longer range.}
\label{picture}
\end{figure}

If $f$ is indefinite, this part needs to be adapted since the region 
$$K_0 \cap \{ (x,y) : f(x,y) \leq \W N' + b\}~,$$ 
which we are summing over, has a different shape.
Let $(x_0,y_0)$ be the point of intersection of 
$\{y=\theta x\}$ and $\{f(x,y) = \W N' + b\}$, and let 
$x_1>0$ be the positive solution to $f(x_1,0)= \W N' + b$. 
Observe that $x_1$ exists since we assumed that $a>0$.
Let $\kappa = 2\max(x_0,x_1)$, and consider the region
$$K_1 = \{(x,y): 0 \leq x \leq \kappa, 
                0 \leq y \leq \min(\theta x, y_0) \}.$$ 
We may now split our given summation as follows:
\begin{align*}
&\frac{1}{N'} \sum_{\substack{x,y \in K_0: \\ f(x,y) \leq \W N' + b}}
  F(g((f(x,y) - b)/\W)\Gamma) \\
&= \frac{1}{N'} \bigg(
  \sum_{x,y \in K_1}
  F(g((f(x,y) - b)/\W)\Gamma) 
- \sum_{\substack{x,y \in K_1 : 
                             \\ f(x,y) > \W N' + b}}
  F(g((f(x,y) - b)/\W)\Gamma)
  \bigg) \\
&= \frac{1}{N'} \sum_{0\leq y \leq y_0}
  \bigg( 
  \sum_{\substack{y/\theta < x \leq \kappa}}
  F(g((f(x,y) - b)/\W)\Gamma)
 -
  \sum_{\substack{y/\theta < x \leq \kappa \\f(x,y) > \W N' + b}}
  F(g((f(x,y) - b)/\W)\Gamma) \bigg)~.
\end{align*}
The next aim is to apply Weyl's inequality for nilsequences to the latter
sums over $x$ for fixed $y$.
Since $x_0 \asymp_D x_1$, the relevant sequences arise as
polynomial
subsequences of $(g(n)\Gamma)_{n \leq A N}$ for some $A=O(1)$.
In \cite{m-definiteforms} we only considered equidistribution properties
of the sequence $(g(n)\Gamma)_{n \leq N}$ and polynomial subsequences
thereof.
For an application of Weyl's inequality, analogous to the one carried out
in \cite[\S18]{m-definiteforms}, it thus remains to justify that for any
$A=O(1)$ there exists a positive integer $A'=O(1)$ such
that $(g(n)\Gamma)_{n \leq A N}$ is totally 
$\delta^{1/A'}$-equidistributed if $(g(n)\Gamma)_{n \leq N}$ is
totally $\delta$-equidistributed.
We shall deduce this assertion from Green and Tao's quantitative Leibman
theorem \cite[Theorem 2.9]{green-tao-polynomialorbits}.

Let $m$ be the dimension of $G/\Gamma$ and let $d$ be the degree of the
filtration $G_{\bullet}$ with respect to which the polynomial sequence
$g$ is defined. Suppose that $(g(n)\Gamma)_{n\leq N}$ is totally
$\delta$-equidistributed, where $\delta = \delta(N) \in (0,1)$ satisfies 
$\delta^{-t} \ll_t N$ for all $t \in \N$. 
By \cite[Proposition 14.3]{m-definiteforms}
there is $B>0$, $B=O_{m,d}(1)$, such that $(\eta \circ g(n))_{n \leq N}$
is totally $\delta^{1/B}$-equidistributed for every horizontal character
$\eta$ of modulus $|\eta| \leq \delta^{-1/B}$.
By \cite[Proposition 14.2(b)]{m-definiteforms} and the definition of the
smoothness norm \cite[Definition 2.7]{green-tao-polynomialorbits}, there
is a positive integer $B'=O_{m,d}(1)$ such that
$$\|k \eta \circ g\|_{C^{\infty}[N]} \geq \delta^{-1/B'}$$
holds for all positive integers $k \leq \delta^{-1/B'}$.
Suppose that $k \eta \circ g: \N_0 \to \R/\Z$ has the
representation $(k \eta \circ g)(n)=\sum_{j=0}^d \alpha_j n^j$.
Then
$$\|k \eta \circ g\|_{C^{\infty}[N]} = \sup_{1 \leq j \leq d} N^j
\|\alpha_j\|
\asymp_A \sup_{1 \leq j \leq d} (AN)^j \|\alpha_j\|
= \| k \eta \circ g \|_{C^{\infty}[AN]}.
$$
Thus 
$\|k \eta \circ g\|_{C^{\infty}[AN]} \gg \delta^{-1/B'}$
for all $k \leq \delta^{-1/B'}$.
By \cite[Theorem 2.9]{green-tao-polynomialorbits}, there therefore is some
integer $B''>0$, $B''=O_{m,d}(1)$ such that $(g(n)\Gamma)_{n\leq AN}$ is
totally $\delta^{1/B''}$-equidistributed.
This completes the proof of the assertion.

\subsection*{Acknowledgement}
The author is grateful to Tim Browning for comments on an earlier draft
of this paper.

\providecommand{\bysame}{\leavevmode\hbox to3em{\hrulefill}\thinspace}

\end{document}